\documentstyle[12s,fl,12pt,twoside]{mybeaut}
\renewcommand{\thesection}{\arabic{section}.}

\newtheorem{conjecture}{Conjecture}

\newtheorem{theorem}{Theorem}

\renewcommand{\refname}{\bf References}
\makeatletter
\renewcommand{\@evenhead}{\thepage \hfil {\leftmark} \hfil}
\renewcommand{\@oddhead}{\hfil {\rightmark} \hfil \thepage}
\makeatother
\input tcilatex

\begin{document}
\count0=1

\title{{\Large \textbf{A Note on Large Cycles in Graphs\\ Around Conjectures of Bondy and Jung}}}
\author{\normalsize  Zhora Nikoghosyan}

\maketitle

\begin{abstract}
Two new sufficient conditions for generalized cycles (including Hamilton and dominating cycles as special cases) in an arbitrary k-connected graph (k=1,2,…) are derived, which prove the truth of Bondy's (1980) famous conjecture for some variants significantly improving the result expected by the given hypothesis. Similarly, two new lower bounds for the circumference (the length of a longest cycle) are established for the reverse hypothesis proposed by Jung (2001).

\noindent {\bf Keywords:} Hamilton cycle, Dominating cycle, Longest cycle, Large cycle.\\
{\bf MSC-class:} 05C38 (primary), 05C45, 05C40 (secondary)
\end{abstract}

\section{Introduction}
We consider only finite undirected graphs without loops or multiple edges. The set of vertices of a graph $G$ is denoted by $V(G)$; the set of edges by $E(G)$. For a subset $S$ of $V(G)$, we denote by $G-S$ the maximum subgraph of $G$ with vertex set $V(G)-S$. For a subgraph $H$ of $G$ we use $G-H$, short for $G-V(H)$. A good reference for any undefined terms is \cite {3}.

     Let $\alpha$ and $\delta$ be the independence number and the minimum degree of a graph $G$, respectively. We define $\sigma_k$ by the  minimum degree sum of any $k$ independent vertices if  $\alpha\ge k$; if $\alpha<k$, we set $\sigma_k=+\infty$. In particular, we have $\sigma_1=\delta$.

     A simple cycle (or just a cycle) of order $t$ (the number of vertices) is a sequence $v_1 v_2...v_t v_1$ of distinct vertices $v_1,...,v_t$ with $v_i v_{i+1}\in E(G)$ for each $i\in\{1,...,t\}$, where $v_{t+1}=v_1$. When $t=1$, the cycle $v_1$ coincides with the vertex $v_1$. So, by this standard definition, all vertices and edges in a graph can be considered as cycles of orders 1 and 2, respectively. Such an extension of the cycle definition allows to avoid unnecessary repetition of a condition  "let $G$ be a graph of order $n\geq3$" in a large number of results. Further, a simple path (or just a path)  of order $t$ is a sequence $v_1 v_2...v_t$ of distinct vertices $v_1,...,v_t$ with $v_i v_{i+1}\in E(G)$ for each $i\in\{1,...,t-1\}$.  

     A graph $G$ is hamiltonian if $G$ contains a Hamilton cycle, i.e., a cycle of order $|V(G)|$.

Now let $Q$ be an arbitrary cycle in $G$. We say that $Q$  is a dominating cycle in $G$  if $V(G-Q)$ is an independent set of vertices.

The first type of generalized cycles, including Hamilton and dominating cycles as special cases, was introduced by Bondy \cite {4}. For a positive integer $\lambda$,  $Q$  is said to be a $D_\lambda$-cycle  if $|H|\leq\lambda-1$ for every component $H$ of $G-Q$. Alternatively, $Q$ is a $D_\lambda$-cycle of $G$ if and only if every connected subgraph of order $\lambda$ of $G$ has at least one vertex with $Q$ in common. Thus a $D_\lambda$-cycle dominates all connected subgraphs of order $\lambda$. By this definition,  $Q$ is a Hamilton cycle if and only if $Q$ is a $D_1$-cycle. Analogously, $Q$ is a dominating cycle if and only if $Q$ is a $D_2$-cycle. 

We now present another two types of more interesting generalized cycles that form the main topic of this paper. For a positive integer $\lambda$, the cycle $Q$  is called a $PD_\lambda$-cycle (PD - Path Dominating) if each path of order at least $\lambda$ in $G$  has at least one  vertex with $Q$ in common.  Similarly, we call the cycle $Q$ a $CD_\lambda$-cycle (CD - Cycle Dominating; introduced in  \cite {13}) if each cycle of order at least $\lambda$ has at least one vertex with $Q$ in common. In fact, a $PD_\lambda$-cycle dominates all paths of order $\lambda$ in $G$; and a $CD_\lambda$-cycle dominates all cycles of order $\lambda$ in $G$. In terms of $PD_\lambda$ and $CD_\lambda$-cycles, $Q$ is a Hamilton cycle if and only if either $Q$ is a $PD_1$-cycle or a $CD_1$-cycle. Further, $Q$ is a dominating cycle if and only if either $Q$ is a $PD_2$-cycle or a $CD_2$-cycle. 

     Throughout the paper, we consider a graph $G$ on $n$ vertices with minimum degree $\delta$ and connectivity $\kappa$. Further, let $C$ be a longest cycle in $G$ with $c=|C|$, and let $\overline{p}$ and $\overline{c}$ denote the orders of a longest path and a longest cycle in $G-C$, respectively. In particular, $C$ is a Hamilton cycle if and only if $\overline{p}\leq 0$ or $\overline{c}\leq 0$. Similarly, $C$ is a dominating cycle if and only if $\overline{p}\leq 1$ or $\overline{c}\leq 1$.

     In 1980, Bondy \cite {4} conjectured a common generalization of some well-known degree-sum conditions for $PD_\lambda$-cycles (called $(\sigma,\overline{p}$)-version) including Hamilton cycles ($PD_1$-cycles) and dominating cycles ($PD_2$-cycles) as special cases.

\begin{conjecture} (Bondy \cite {4},1980): $(\sigma, \overline{p})$-version \\
Let $C$ be a longest cycle in a $\lambda$-connected $(1\leq\lambda\leq \delta)$ graph $G$ of order $n$. If $\sigma_{\lambda+1}\geq n+\lambda(\lambda-1)$, then $\overline{p}\leq \lambda-1$.
\end{conjecture}

     Parts of Conjecture 1 were proved by Ore  \cite {15} $(\lambda=1)$, Bondy [4] $(\lambda=2)$ and Zou [17] $(\lambda=3)$. For the general case, Conjecture 1 is still open.

     The long cycles analogue (so called reverse version) of Bondy’s conjecture (Conjecture 1) can be formulated as follows. 

\begin{conjecture}: (reverse, $\sigma$, $\overline{p})$-version \\
Let $C$ be a longest cycle in a $\lambda$-connected $(1\leq\lambda\leq\delta)$ graph $G$. If $\overline{p}\geq \lambda-1$, then $c\geq \sigma_\lambda-\lambda(\lambda-2)$.
\end{conjecture}

     Parts of Conjecture 2 were proved by Dirac  \cite {6} $(\lambda=1)$, Bondy  \cite {2}, Bermond  \cite {1}, Linial  \cite {11} $(\lambda=2)$, Fraisse, Yung  \cite {8} $(\lambda=3)$ and Chiba, Tsugaki, Yamashita  \cite {5} $(\lambda=4)$.

     Note that the initial motivations of Conjecture 1 and Conjecture 2 come from their minimal degree versions - the most popular and much studied versions, which also remain unsolved.  

\begin{conjecture} (Bondy \cite {4},1980): $(\delta, \overline{p})$-version \\
Let $C$ be a longest cycle in a $\lambda$-connected $(1\leq\lambda\leq \delta)$ graph $G$ of order $n$. If $\delta\geq\frac{n+2}{\lambda+1}+\lambda-2$, then $\overline{p}\leq \lambda-1$.
\end{conjecture}

\begin{conjecture}: (Jung \cite{10}, 2001): (reverse, $\delta$, $\overline{p})$-version \\
Let $C$ be a longest cycle in a $\lambda$-connected $(1\leq\lambda\leq\delta)$ graph $G$. If $\overline{p}\geq \lambda-1$, then $c\geq\lambda(\delta-\lambda+2)$.
\end{conjecture}

Parts of Conjecture 3 were proved for $\lambda=1,2,3$. 

\ \

\noindent $(a)\ \ \kappa\geq1, \ \ \delta\geq\frac{n}{2} \ \ \ \ \ \ \  \Longrightarrow \ \ \ \  \overline{p}\leq 0  \ \ \ (Dirac \cite{6},1952), $ 

\noindent $(b)\ \ \kappa\geq2, \ \ \delta\geq\frac{n+2}{3} \ \ \ \ \  \Longrightarrow \ \ \ \  \overline{p}\leq 1  \ \ \ (Nash-Williams \cite{12},1971), $ 

\noindent$(c)\ \ \kappa\geq3, \ \ \delta\geq\frac{n+6}{4} \ \ \ \ \  \Longrightarrow \ \ \ \  \overline{p}\leq 2  \ \ \ (Fan \cite{7},1987). $ 

\ \

Parts of Conjecture 4 were proved for $\lambda=1,2,3,4$.

\ \

\noindent $(d) \ \ \kappa\geq 1, \ \  \overline{p}\geq 0 \ \ \ \ \ \ \  \Longrightarrow \ \ \ \   c\geq \delta+1 \  \ \ \ (Dirac \cite{6},1952), $

\noindent $(e) \ \ \kappa\geq 2, \ \  \overline{p}\geq 1 \ \  \ \ \ \ \  \Longrightarrow \ \ \ \   c\geq 2\delta \ \ \ \ \ \ \ \ (Dirac \cite{6},1952),$
 
\noindent $(f) \ \ \kappa\geq 3, \ \  \overline{p}\geq 2 \ \ \ \ \ \ \  \Longrightarrow \ \ \ \   c\geq 3\delta-3 \ \ \ (Voss, Zuluaga  \cite{16},1977),$
 
\noindent $(g) \ \ \kappa\geq 4, \ \  \overline{p}\geq 3 \ \ \ \ \ \ \  \Longrightarrow \ \ \ \   c\geq 4\delta-8 \ \ \ (Jung \cite{9},1990).$

\ \

As for the geneal case, parts of Conjecture 1 were proved for $\lambda=1,2,3$.

\ \

\noindent $(h) \ \ \kappa\geq 1, \ \ \sigma_2\geq n \ \ \ \ \ \ \ \ \ \ \Longrightarrow \ \ \ \  \overline{p}\leq 0 \ \ \ \ (Ore \cite{15}, 1960), $

\noindent $(i) \ \ \kappa\geq 2, \ \ \sigma_3\geq {n+2} \ \ \ \ \ \ \Longrightarrow \ \ \ \  \overline{p}\leq 1 \ \ \ \ (Bondy \cite{4}, 1980), $ 

\noindent $(j) \ \ \kappa\geq 3, \ \ \sigma_4\geq {n+6} \ \ \ \ \ \ \Longrightarrow \ \ \ \  \overline{p}\leq 2 \ \ \ \ (Zou \cite{17}, 1987). $

\ \

Finally, parts of Conjecture 2 were proved for $\lambda=1,2,3,4$.

\ \

\noindent $(k) \ \ \kappa\geq 1, \ \  \overline{p}\geq 0 \ \ \ \ \  \Longrightarrow \ \ \ \   c\geq \sigma_1+1 \  \ \ \ (Dirac \cite{6},1952),$

\noindent $(l) \ \ \kappa\geq 2, \ \  \overline{p}\geq 1 \ \ \ \ \ \  \Longrightarrow \ \ \ \   c\geq \sigma_2 \  \ \ \ (Bondy \cite{2}, 1971; Bermond \cite{1}, 1976; Linial \cite{11}, 1976)$,

\noindent $(m) \ \ \kappa\geq 3, \ \  \overline{p}\geq 2 \ \ \ \  \Longrightarrow \ \ \ \   c\geq \sigma_3-3 \  \ \ \ (Fraisse, Jung \cite{8},1989),$

\noindent $(n) \ \ \kappa\geq 4, \ \  \overline{p}\geq 3 \ \ \ \ \  \Longrightarrow \ \ \ \   c\geq \sigma_4-8 \  \ \ \ (Chiba, Tsugaki, Yamashita \cite{5}, 2014)$.

\ \

Note that $CD_\lambda$-cycles are more suitable for research than $PD_\lambda$-cycles since cycles in $G-C$ are more symmetrical than paths in view of the connections between $G-C$ and $CD_\lambda$-cycles. This is the main reason why some minimum degree versions of Conjectures 1 and 2 have been solved just for $CD_\lambda$-cycles. 

     According to above arguments, it is natural to consider the exact analogues of Bondy’s generalized conjecture (Conjecture 1) and its reverse version (Conjecture 2) for $CD_\lambda$-cycles which we call $(\sigma,\overline{c})$ and $(reverse, \sigma, \overline{c})$-versions, respectively.

\begin{conjecture}: $(\sigma, \overline{c})$-version \\
Let $C$ be a longest cycle in a $\lambda$-connected $(1\le\lambda\le\delta)$ graph $G$ of order $n$. If $\sigma_{\lambda+1}\geq n+\lambda(\lambda-1)$, then $\overline{c}\leq \lambda-1$.
\end{conjecture}

\begin{conjecture}: (reverse, $\sigma$, $\overline{c})$-version \\
Let $C$ be a longest cycle in a $\lambda$-connected $(1\le\lambda\le\delta)$ graph. If $\overline{c}\geq \lambda-1$, then $c\geq \sigma_\lambda-\lambda(\lambda-2)$.
\end{conjecture}

   In 2009, the author proved  \cite {14} the validity of minimum degree versions of Conjectures 5 and 6.

\begin{theorem} (\cite {14}, 2009): $(\delta, \overline{c})$-version \\
Let $C$ be a longest cycle in a $\lambda$-connected $(1\le\lambda\le \delta)$ graph $G$ of order $n$. If $\delta\geq \frac{n+2}{\lambda+1}+\lambda-2$, then $\overline{c}\leq \lambda-1$.
\end{theorem}

\begin{theorem} (\cite {14}, 2009): (reverse, $\delta$, $\overline{c})$-version \\
Let $C$ be a longest cycle in a $\lambda$-connected $(1\le\lambda\le\delta)$ graph. If $\overline{c}\geq \lambda-1$, then $c\geq \lambda(\delta-\lambda+2)$.
\end{theorem}

    Actually, in \cite{14} it was proved a significantly stronger result than Theorem 1 by showing that the conclusion $\overline{c}\le\lambda-1$ in Theorem 1 can be strengthened to $\overline{c}\le\min\{\lambda-1,\delta-\lambda\}$, called $\overline{c}$-improvement.

\begin{theorem} (\cite {14}, 2009): $(\delta, \overline{c})$-version, $\overline{c}$-improvement \\
Let $C$ be a longest cycle in a $\lambda$-connected $(1\le\lambda\le\delta)$ graph $G$ of order $n$. If $\delta\geq \frac{n+2}{\lambda+1}+\lambda-2$, then $\overline{c}\leq \min\{ \lambda-1,\delta-\lambda\}$.
\end{theorem}

 Analogously, the condition $\overline{c}\geq \lambda-1$ in Theorem 2 was weakened  \cite {14} to $\overline{c}\geq \min\{\lambda-1,\delta-\lambda+1\}$. 

\begin{theorem} (\cite {14}, 2009): (reverse, $\delta$, $\overline{c})$-version, $\overline{c}$-improvement \\
Let $C$ be a longest cycle in a $\lambda$-connected $(1\le\lambda\le\delta)$ graph $G$. If $\overline{c}\geq\min\{ \lambda-1,\delta-\lambda+1\}$, then $c\geq \lambda(\delta-\lambda+2)$.
\end{theorem}

In this paper we present four new analogous improvements  of Theorems 1,2,3,4 inspiring a number of new conjectures in forms of improvements of initial generalized conjectures of Bondy and Jung.

\section{Results}

First, we prove that the connectivity condition $\kappa\ge\lambda$ in Theorem 1 can be weakened to $\kappa\ge\min\{\lambda,\delta-\lambda+1\}$.

\begin{theorem}: $(\delta, \overline{c})$-version, $ \kappa$-improvement \\
Let $C$ be a longest cycle in a graph $G$ of order $n$ and $\lambda$ a positive integer with $1\le\lambda\le\delta$. If $\kappa\geq \min\{\lambda,\delta-\lambda+1\}$ and $\delta\geq \frac{n+2}{\lambda+1}+\lambda-2$, then $\overline{c}\leq \lambda-1$.
\end{theorem}

Analogously, we prove that the connectivity condition $\kappa\ge\lambda$  in Theorem 2 can be weakened to  $\kappa\ge \min\{\lambda,\delta-\lambda+2\}$.

\begin{theorem}: (reverse, $\delta$, $\overline{c})$-version, $\kappa$-improvement \\
Let $C$ be a longest cycle in a graph $G$ and $\lambda$ a positive integer with  $1\le\lambda\le\delta$. If $\kappa\ge\min\{\lambda, \delta-\lambda+2\}$ and $\overline{c}\ge \lambda-1$, then $c\ge\lambda(\delta-\lambda+2)$.  
\end{theorem}

Next, we prove that the conclusion $\overline{c}\le\lambda-1$ in Theorem 5 can be strengthened to $\overline{c}\le\min\{\lambda-1,\delta-\lambda\}$.

\begin{theorem}: $(\delta, \overline{c})$-version, $(\overline{c}, \kappa$)-improvement \\
Let $C$ be a longest cycle in a graph $G$ of order $n$ and $\lambda$ a positive integer with $1\le\lambda\le\delta$. If $\kappa\geq \min\{\lambda,\delta-\lambda+1\}$ and $\delta\geq \frac{n+2}{\lambda+1}+\lambda-2$, then $\overline{c}\le\min\{\lambda-1,\delta-\lambda\}$.
\end{theorem}

Finally, we prove that the condition $\overline{c}\ge \lambda-1$ in Theorem 6 can be weakened to $\overline{c}\ge \min\{\lambda-1,\delta-\lambda+1\}$.

\begin{theorem}: (reverse, $\delta$, $\overline{c})$-version, $(\overline{c},\kappa)$-improvement \\
Let $C$ be a longest cycle in a graph $G$ and $\lambda$ a positive integer with  $1\le\lambda\le\delta$. If $\kappa\ge\min\{\lambda, \delta-\lambda+2\}$ and $\overline{c}\ge \min\{\lambda-1,\delta-\lambda+1\}$, then $c\ge \lambda(\delta-\lambda+2)$.  
\end{theorem}

\section{Generalized Improvements of Conjectures of Bondy and Jung }

Motivated by Theorems 5,6,7,8 (minimum degree versions) with Conjectures 1 and 2, in this section we propose their exact analogs in terms of degree sums as generalized improvements of Conjectures of Bondy and Jung.

\begin{conjecture}: $(\sigma, \overline{c})$-version, $(\overline{c}, \kappa$)-improvement \\
Let $C$ be a longest cycle in a graph $G$ of order $n$ and $\lambda$ a positive integer. If $\kappa\geq \min\{\lambda,\delta-\lambda+1\}$ and $\sigma_{\lambda+1}\geq n+\lambda(\lambda-1)$, then $\overline{c}\leq \min\{\lambda-1,\delta-\lambda\}$.
\end{conjecture}

\begin{conjecture}: (reverse, $\sigma$, $\overline{c})$-version, $(\overline{c},\kappa$)-improvement \\
Let $C$ be a longest cycle in a graph $G$ and $\lambda$ a positive integer. If $\kappa\ge\min\{\lambda, \delta-\lambda+2\}$ and $\overline{c}\ge \min\{\lambda-1,\delta-\lambda+1\}$, then $c\ge \sigma_\lambda-\lambda(\lambda-2)$.  
\end{conjecture}

\begin{conjecture}: $(\sigma, \overline{p})$-version, $(\overline{p},\kappa$)-improvement \\
Let $C$ be a longest cycle in a graph $G$ of order $n$ and $\lambda$ a positive integer. If $\kappa\geq \min\{\lambda,\delta-\lambda+1\}$ and $\sigma_{\lambda+1}\geq n+\lambda(\lambda-1)$, then $\overline{p}\le\min\{\lambda-1,\delta-\lambda\}$.
\end{conjecture}

\begin{conjecture}: (reverse, $\sigma$, $\overline{p})$-version, $(\overline{p},\kappa)$-improvement \\
Let $C$ be a longest cycle in a graph $G$ and $\lambda$ a positive integer. If $\kappa\ge\min\{\lambda, \delta-\lambda+2\}$ and $\overline{p}\ge \min\{\lambda-1,\delta-\lambda+1\}$, then $c\ge \sigma_\lambda-\lambda(\lambda-2)$.
\end{conjecture}

 \section{Proofs}

\noindent\textbf{Proof of Theorem 7}.   We shall prove that $\overline{c}\leq \min\{\lambda-1,\delta-\lambda\}$ under the conditions
$$
\kappa\geq \min\{\lambda,\delta-\lambda+1\},\ \ \delta\geq\frac{n+2}{\lambda+1}+\lambda-2
$$
for each $1\le\lambda\le \delta$. If  $\min\{\lambda,\delta-\lambda+1\}=\lambda$, that is $\lambda\leq \lfloor \frac{\delta+1}{2}\rfloor$, then we have $\kappa\geq \lambda$ and the conclusion $\overline{c}\leq \lambda-1$ in Theorem 7 follows from Theorem 1 for all $\lambda=1,2,...,\lfloor\frac{\delta+1}{2}\rfloor$. Now let $\min\{\lambda,\delta-\lambda+1\}=\delta-\lambda+1$, that is $\lambda\geq \lfloor \frac{\delta+2}{2}\rfloor$. To conclude the proof, it remains to show that 
\begin{equation}
\kappa\geq\delta-\lambda+1,\  \delta\geq\frac{n+2}{\lambda+1}+\lambda-2\ \  \Rightarrow \ \ \overline{c}\leq \delta-\lambda \ \ \  \left(\lambda=\delta,\delta-1,..., \left\lfloor \frac{\delta+2}{2}\right\rfloor\right).
\end{equation}

 Put $\delta-\lambda+1=\mu$. Acording to this notation, (1) is equivalent to 
\begin{equation}
\kappa\geq\mu,\  \delta\geq\frac{n+2}{\delta-\mu+2}+\delta-\mu-1 \ \  \Rightarrow \ \ \overline{c}\leq \mu-1 \ \ \  \left(\mu=1,2,..., \left\lfloor \frac{\delta+1}{2}\right\rfloor\right).
\end{equation}

In (2),  the inequality
$$
\delta\geq\frac{n+2}{\delta-\mu+2}+\delta-\mu-1
$$
is equivalent to
$$
\delta\geq\frac{n+2}{\mu+1}+\mu-2,
$$
implying that (2) is equivalent to

\begin{equation}
\kappa\geq\mu,\  \delta\geq\frac{n+2}{\mu+1}+\mu-2 \ \  \Rightarrow \ \ \overline{c}\leq \mu-1 \ \ \  \left(\mu=1,2,..., \left\lfloor \frac{\delta+1}{2}\right\rfloor\right). 
\end{equation}

Observing that (3) follows from Theorem 1 immediately, we obtain
$$
(1)\equiv (2)\equiv (3)\ \Leftarrow \ "Theorem \ 1".
$$

Theorem 7 is proved.\ \ \ \  \endproof   \\

\noindent\textbf{Proof of Theorem 5}.  Let $G$ be a graph with
$$
\kappa\ge\min\{\lambda,\delta-\lambda+1\}, \  \  \delta\ge\frac{n+2}{\lambda+1}+\lambda-2
$$
for each $1\le\lambda\le\delta$. We shall prove that $\overline{c}\le\lambda-1$. If 
 $\min\{\lambda,\delta-\lambda+1\}=\lambda$, that is 
$\lambda\leq \lfloor \frac{\delta+1}{2}\rfloor$,
then we have
$\kappa\ge\lambda$ and the conclusion $\overline{c}\le\lambda-1$  in Theorem 5 follows from Theorem 1 for all 
$\lambda=1,2,...,\lfloor \frac{\delta+1}{2}\rfloor$. Now let $\min\{\lambda,\delta-\lambda+1\}=\delta-\lambda+1$,
 that is $\delta-\lambda\le\lambda-1$ and $\lambda\geq \lfloor \frac{\delta+2}{2}\rfloor$. To conclude the proof, it remains to show that
\begin{equation}
\kappa\ge\delta-\lambda+1, \ \ \delta\ge\frac{n+2}{\lambda+1}+\lambda-2  \ \  \Rightarrow \ \ \overline{c}\le \lambda-1 \ \  \left(\lambda=\delta,\delta-1,..., \left\lfloor \frac{\delta+2}{2}\right\rfloor\right).
\end{equation}
By Theorem 7,
\begin{equation}
\kappa\ge\delta-\lambda+1, \ \ \delta\ge\frac{n+2}{\lambda+1}+\lambda-2  \ \  \Rightarrow \ \ \overline{c}\le \delta-\lambda \ \  \left(\lambda=\delta,\delta-1,..., \left\lfloor \frac{\delta+2}{2}\right\rfloor\right).
\end{equation}
Recalling that $\delta-\lambda\le\lambda-1$, we can weaken the conclusion $\overline{c}\le \delta-\lambda$ in (5) to $\overline{c}\le \lambda-1$. Then (4) holds immediately, which completes the proof of Theorem 5.   \ \ \ \  \endproof   \\
  
\noindent\textbf{Proof of Theorem 8}.  Let $G$ be a graph with
$$
\kappa\geq \min\{\lambda,\delta-\lambda+2\}, \ \ \overline{c}\geq \min\{\lambda-1,\delta-\lambda+1\}
$$
for each $1\leq \lambda\leq \delta$. We shall prove that  $c\geq\lambda(\delta-\lambda+2)$. If $\lambda=1$, then the result follows from the fact that each graph has a cycle of length at least $\delta+1$ \cite{6}. Let $\lambda\geq 2$. Further, if  $\min\{\lambda,\delta-\lambda+2\}=\lambda$, then we are done by Theorem 2. Now let $\min\{\lambda,\delta-\lambda+2\}=\delta-\lambda+2$, that is  $\lambda\geq \lfloor \frac{\delta+3}{2}\rfloor$. Then it remains to prove that 
\begin{equation}
\kappa\geq \delta-\lambda+2,\ \ \overline{c}\geq\delta-\lambda+1 \ \Rightarrow \  c\geq\lambda(\delta-\lambda+2) \ \ \  \left(\lambda=\delta,\delta-1,..., \left\lfloor \frac{\delta+3}{2}\right\rfloor\right).
\end{equation}

Put $\delta-\lambda+2=\mu$. By this notation,  the statement (6) is equivalent to  

\begin{equation}
\kappa\geq \mu \ \ \overline{c}\geq\mu-1 \ \Rightarrow \  c\geq\mu(\delta-\mu+2) \ \ \  \left(\mu=2,3,..., \left\lfloor \frac{\delta+2}{2}\right\rfloor\right),
\end{equation}
which follows from Theorem 2 immediately. So, $(6)\equiv (7) \Leftarrow "Theorem \ 2"$.   Theorem 8 is proved.      \ \ \endproof \\

 \noindent\textbf{Proof of Theorem 6}.  Let $G$ be a graph with $\kappa\ge\min\{\lambda,\delta-\lambda+2\}$ and $\overline{c}\ge\lambda-1$ for each $1\le\lambda\le \delta$. We shall prove that $c\ge \lambda(\delta-\lambda+2)$. If $\lambda=1$, then the result follows from the fact that in each graph, $c\ge\delta+1$.
 Let $\lambda\ge 2$. Further, if $\min\{\lambda,\delta-\lambda+2\}=\lambda$, then we are done by Theorem 2, Now let $\min\{\lambda,\delta-\lambda+2\}=\delta-\lambda+2$,
 that is $\lambda\ge \delta-\lambda+2$ and
 $\lambda\ge \lfloor \frac{\delta+3}{2}\rfloor$. It remains to prove that 
\begin{equation}
\kappa\geq \delta-\lambda+2, \ \ \overline{c}\geq\lambda-1 \ \Rightarrow \  c\geq\lambda(\delta-\lambda+2) \ \ \  \left(\lambda=\delta,\delta-1,..., \left\lfloor \frac{\delta+3}{2}\right\rfloor\right).
\end{equation}
By Theorem 8, we have
\begin{equation}
\kappa\geq \delta-\lambda+2, \ \ \overline{c}\geq\delta-\lambda+1 \ \Rightarrow \  c\geq\lambda(\delta-\lambda+2) \ \ \  \left(\lambda=\delta,\delta-1,..., \left\lfloor \frac{\delta+3}{2}\right\rfloor\right).
\end{equation}
Recalling that $\lambda\ge\delta-\lambda+2$, we can strengthen the condition $\overline{c}\ge\delta-\lambda+1$ in (9) to $\overline{c}\ge\lambda-1$. This completes the proof of (8) and Theorem 6.  \ \ \endproof.

\noindent  Institute for Informatics and Automation Problems of NAS RA\\
P. Sevak 1, Yerevan 0014, Armenia\\
e-mail: zhora@iiap.sci.am

\end{document}